\documentclass[12pt]{amsart}
\usepackage{amsbsy,amsfonts,amssymb,amsmath,graphicx}
\begin{document}

\newtheorem{theorem}{Theorem}
\newtheorem{proposition}{Proposition}
\newtheorem{lemma}{Lemma}
\newtheorem{corollary}{Corollary}
\newtheorem{definition}{Definition}
\newtheorem{remark}{Remark}


\numberwithin{equation}{section}
\numberwithin{theorem}{section}
\numberwithin{proposition}{section}
\numberwithin{lemma}{section}
\numberwithin{corollary}{section}
\numberwithin{definition}{section}
\numberwithin{remark}{section}

\newcommand{\beq}{\begin{equation}}
\newcommand{\eeq}{\end{equation}}


\newcommand{\ren}{{\mathbb R}^N}
\newcommand{\re}{{\mathbb R}}
\newcommand{\n}{\nabla}
\newcommand{\iy}{\infty}
\newcommand{\pa}{\partial}
\newcommand{\fp}{\noindent}
\newcommand{\ms}{\medskip\vskip-.1cm}
\newcommand{\mpb}{\medskip}
\newcommand{\ssk}{\smallskip}


\renewcommand{\a}{\alpha}
\renewcommand{\b}{\beta}
\newcommand{\g}{\gamma}
\newcommand{\G}{\Gamma}
\renewcommand{\d}{\delta}
\newcommand{\D}{\Delta}
\newcommand{\e}{\varepsilon}
\newcommand{\vp}{\varphi}
\renewcommand{\l}{\lambda}
\renewcommand{\o}{\omega}
\renewcommand{\O}{\Omega}
\newcommand{\s}{\sigma}
\renewcommand{\t}{\tau}
\renewcommand{\th}{\theta}
\newcommand{\z}{z}
\newcommand{\wx}{\widetilde x}
\newcommand{\wt}{\widetilde t}
\newcommand{\noi}{\noindent}
\newcommand{\AAA}{{\mathbf   A}}
\newcommand{\BB}{{\mathbf  B}}
\newcommand{\CC}{{\mathbf  C}}

\title[Positivity]
{On convergence in smooth gradient \\
 systems with branching of 
 equilibria}

\author
{V.A. Galaktionov, S.I. Pohozaev,
 and A.E. Shishkov}

\address{Department of Math. Sci., University of Bath,
 Bath BA2 7AY, UK}
\email{vag@maths.bath.ac.uk}

\address{Steklov Mathematical Institute,
 Gubkina St. 8, 119991 Moscow, RUSSIA}
\email{pokhozhaev@mi.ras.ru}

\address{Institute of Applied Mathematics and Mechanics
of NAS of Ukraine, R. Luxemburg str. 74, 83114 Donetsk, UKRAINE}
\email{shishkov@iamm.ac.donetsk.ua}


\date{\today}

\begin{abstract}
Our basic model is a semilinear elliptic equation with a coercive
$C^1$ nonlinearity, $\D \psi+f(\psi)=0$ in $\O$, $\psi=0$ on
$\partial \O$, where $\O \subset \ren$ is a bounded smooth domain.
Our main hypothesis ${\bf (H_R)}$ on the {\em resonance branching}
 is as follows: if a branching of equilibria  occurs at any point
$\psi$ with a $k$-dimensional kernel of the linearized  operator
$\D + f'(\psi)I$, the branching subset $S_k$ at $\psi$ is a
locally smooth  $k$-dimensional manifold.

 Using the corresponding
parabolic flow
 $$
u_t= \D u + f(u),
 $$
 with bounded initial data $u_0$,
 we then prove that, under the (${\bf H_R}$),
    the subset of equilibria $\Phi=\{\psi\}$ is evolutionary complete,
 i.e.,  any evolution orbit converges as $t \to \infty$
to a single element of $\Phi$. We also treat the non-coercive case
$f(u)=|u|^{p-1}u$ for $p \in (1, \frac{N+2}{N-2})$, where the
elliptic problem is known to admit at least a countable subset of
solutions $\{\psi_k\}$ due to the classical Lusternik-Schnirel'man
category theory. The results are extended to higher-order elliptic
operators $-(-\D)^m u + |u|^{p-1}u$ with Dirichlet conditions on
$\partial \O$, $1<p<\frac{N+2m}{N-2m}$, $m \ge 2$, with the
corresponding adaptation of the parabolic flows.

For $N=1$, the first result on stabilization to a single
equilibrium is due to T.I.~Zelenyak (1968). We show that
Zelenyak's approach based on Lyapunov's function analysis can be
extended to general gradient systems in Hilbert spaces with a
similar smooth resonance branching. We also cover the case of
their asymptotically small non-autonomous perturbations. In
general, the developed approach represents an alternative (and
improved) method on stabilization in Hale--Raugel (1992) and other
later similar techniques in gradient system theory.


\end{abstract}


\maketitle


\begin{center}
{\em Dedicated to  the memory of Professor T.I.~Zelenyak}
\end{center}

\ssk\ssk\ssk\ssk

This is an extended version of the paper \cite{GPSZel}, containing extra comments and updated references.

\section{\bf Introduction: the convergence problem and evolution completeness}
\label{Sect1}

\subsection{Semilinear elliptic
equations and the convergence problem}

Let  $\O \subset \ren$ be a bounded domain with the smooth
boundary $\partial \O$. Our basic model is the classical semilinear
coercive elliptic equation from combustion theory
 \beq
 \label{Par.11}
V'(\psi) \equiv  \D \psi +
 f(\psi) 
 =0 \,\,\, \mbox{in} \,\,\, \O, \,\,\, \psi=0
 \,\,\,
 \mbox{on} \,\,\, \partial \O,
 \eeq
where $f$ is a given function satisfying
 \beq
 \label{f.cond}
 f \in C^1(\re), \quad f'(u) < 0 \,\,\, \mbox{for \,\,$|u| \gg
 1$}.
 \eeq
This  non-optimal coercivity condition is sufficient for demonstrating
the approach.
 For instance, $f(u)=u - u^3$ (an analytic nonlinearity) or $f(u)=|u|^{p-1}u- |u|^{q-1}u$
 with any $q>p>1$. 
 The operator in (\ref{Par.11}) is potential and coercive, so
 the classical variational theory establishes existence of
 solutions as critical points of the corresponding functional
 in $H^1_0(\O)$,
 \beq
 \label{FF.1}
 \mbox{$
 V(\psi) = - \frac 12 \, \int |D \psi|^2 \,{\mathrm d} x + \int F(\psi)\,{\mathrm d}
 x, $}
  \eeq
  where $F'=f$; see e.g., \cite{KrasZ}. If $f$ does not satisfy the
  coercivity condition (\ref{f.cond}), then, in the most well-known 
  power case
  \beq
  \label{f.f.1}
  \mbox{$
  f(u)=|u|^{p-1}u \quad \mbox{with \, $1<p< p_S= \frac{N+2}{N-2}$},
  $}
  \eeq
 versions of the classical Lusternik--Schnirel'man (L-S) theory \cite{Lus0}--\cite{LusS1}
  (see e.g., \cite{Clark} and
\cite{Poh0, PohFM}) establish existence of at least a countable
subset of critical points of the corresponding functional
(\ref{FF.1}).

 By $\Phi= \{\psi\}$ we denote  the subset of all bounded
solutions of the elliptic problem (\ref{Par.11}) in both coercive
and non-coercive  cases. The {\em evolution completeness} of the
subset $\Phi$ of solutions of (\ref{Par.11}) is defined by
introducing the corresponding parabolic equation
 \beq
 \label{Par.12}
u_t =V'(u) \equiv  \D u + f(u)  \,\,\,\, \mbox{in} \,\,\, \O
\times \re_+, \,\,\, u=0
 \,\,\,
 \mbox{on} \,\,\, \partial \O \times \re_+,
 \eeq
 with arbitrary bounded initial data $u_0$.
It follows from the classical parabolic theory \cite{Fr}, that
under hypotheses (\ref{f.cond}) the orbit $\{u(\cdot,t), t>0\}$ is
uniformly bounded. For such smooth gradient systems,  one can
define its $\o$-limit set $\o(u_0)$, which is non-empty, compact,
invariant and connected in the topology of $L^2(\O)$ (or $C(\O)$),
and
 \beq
 \label{Om1}
 \o(u_0) \subseteq \Phi;
 \eeq
 see \cite{Ha} and \cite[pp.~483-487]{Sell}.
The notion of the evolution completeness of the subset  $\Phi$ of
nonlinear eigenfunctions was introduced in \cite{GalC} for the
porous medium equation
 $$
 \mbox{$
 u_\t = \D (|u|^{m-1}u) + \frac 1{m-1} \, u \quad (m>1)
 $}
 $$
 in a bounded domain. In this case,
$0$ is not included into $\Phi$ (as usual, 0 is not an
eigenfunction) and this determines some specific difficulties of
the asymptotic analysis.

In the present problem, $0 \in \Phi$, so that
the evolution completeness of $\Phi$ just means that, for the parabolic flow (\ref{Par.12}),
 \beq
\label{evC}
\mbox{for any data $u_0$, \, $\exists$ a $ \psi \in \Phi$ such that   $\o(u_0)= \{\psi\}$}.
 \eeq
Indeed, this is the classical {\em convergence  problem} in the
theory of dynamical systems. Obviously, (\ref{evC}) is true for
any smooth gradient system, if $\Phi$ consists of isolated points
(a straightforward consequence of the connectedness of $\o(u_0)$),
i.e.,
 \beq
\label{bb1}
\Phi \,\,\, \mbox{is discrete} \,\, \Longrightarrow
 \,\, \mbox{evolution completeness (\ref{evC}) holds}.
 \eeq


\ssk

\noi{\bf On branching situations: two approaches.} The case, where
$\Phi$ is not discrete and consists of continuous families due to
branching of equilibrium points was for a long time a problem of
concern. Existence of various branches of stationary solutions is
established by classical variational techniques; see \cite{KrasZ}
and \cite{Kiel} for  applications to semilinear elliptic problems.
In some cases, the fibering method
\cite{Poh0, PohFM} shows a clear picture of the connection
of  the number of branches with the algebraic non-monotonicity of
$f(u)$, so,  in general, $\Phi$ may contain continuous sub-families;
see  examples in Section \ref{SectRes}.

The first approach to convergence in parabolic problems in the
presence of branching  was formulated in 1968 by {\em
T.I.~Zelenyak} \cite{Zel}, who proved (\ref{evC}) for arbitrary
one-dimensional uniformly parabolic equations with smooth
coefficients
 \beq
\label{Par1}
u_t = a(x,u,u_x) u_{xx} + b(x,u,u_x)
 \eeq
with arbitrary nonlinear boundary
conditions.
 Ten years later, in 1978, a
similar result  was proved by {\em H.~Matano} \cite{Ma0} by
Sturm's zero set argument, which turned out to be applied to
arbitrary smooth fully nonlinear uniformly parabolic equations for
$N=1$
 \beq
\label{Par2} u_t = F(x,u,u_x,u_{xx}).
 \eeq
Matano's approach is principally one-dimensional, since Sturm's
Theorem on zero sets is not available in higher dimensions,
 cannot be extended to parabolic PDEs in $\ren$ for $N \ge 2$.

Another direction of the  development of the stability theory of
gradient systems is associated with  the {\em
\L{ojasiewicz}--Simon} inequality and technique. It is well-known
that, for gradient dynamical systems in a real Banach space $E$
 \beq
 \label{DS.1}
 u_t = V'(u) \quad \mbox{for} \,\,\, t>0, \quad u(0)=u_0 \in E,
 \eeq
 where $V$ is a $C^2$ functional, with the $C^1$ Frechet
 derivative $V':E \to E'$, $E'$ being the dual space,
 the
\L{ojasiewicz}--Simon inequality is an effective tool of studying
of the stabilization phenomena. In particular, it completely
settles the case of analytic nonlinearities; see references in
\cite{Chill}.
Given an equilibrium point $\psi$, with $V'(\psi)=0$, this classical
inequality has the form
 \beq
 \label{LS1}
 |V(v)-V(\psi)|^{1-\theta} \le c \| V'(v)\|_{E'}
 \eeq
 in a neighbourhood of $\psi$,
 where $c>0$ is a constant and $\theta \in (0, \frac 12]$
is called the {\em \L{ojasiewicz} exponent}; see basic references and
historical comments in \cite{Chill}. If (\ref{LS1}) with a fixed
exponent $\theta$ holds for the {\rm critical manifold} \cite[Thm.
3.10]{Chill}, then, under natural assumptions on the spectral
properties of linearized operators (all those hold for our
elliptic case), it is possible to
 establish:

\ssk
 \noi (i) that the $\o$-limit set $\o(u_0)$  of any
 bounded orbit consists of a unique equilibrium, i.e.,
(\ref{evC}) holds, and

\ssk

 \noi (ii) the rate of convergence to $\psi$ as $t \to
\infty$.

\ssk

\noi For instance,  in the  case of finite
regularity, the best possible constant $\theta = \frac 12$ guarantees
 the {\em exponential} convergence,  \cite[Thm. 1.1]{Haraux},
 \beq
 \label{th2}
 \|u(t)-\psi\|_E \le c \, {\mathrm e}^{- \d t} \quad \mbox{as} \,\,\,
 t \to \infty \quad (\d>0).
  \eeq
 In general, roughly speaking (see examples in \cite{Chill})
  \beq
  \label{th1}
  \mbox{$
  \theta \in (0, \frac 12) \,\, \Longrightarrow \,\,
 \|u(t)-\psi\|_E \le c \, t^{- \frac \theta{1-2 \theta}} \quad \mbox{as} \,\,\,
 t \to \infty.
 $}
  \eeq

Concerning applications of the \L{ojasiewicz}--Simon approach to
parabolic equations like (\ref{Par.12}), it does
not assume any restrictions of the  dimension $N$, though, in presence of
branching,
 checking the validity of the corresponding inequality on the singular subset becomes an extremely
difficult problem; see examples in \cite{Chill}. On the other hand, in Section \ref{SectRate}
we show that the \L{ojasiewicz}--Simon inequality (\ref{LS1})  does
not apply to general $C^\infty$ dynamical systems,
i.e., the exponent $\theta \in(0,\frac 12]$ does not always exist.
Actually, this means that
this kind of  analysis needs other non-rational moduli of continuity
in the inequality (\ref{LS1}).


\ssk

\noi{\bf Counterexamples and the resonance branching hypothesis.}
It is well-known for a long time (see \cite{Palis}) that, for
general non-analytic gradient systems on the plane $\re^2$, the
convergence result (\ref{evC}) is not true. The idea of such a
construction
  has been
extended to prove existence of non-convergent orbits for special
 parabolic equations like (\ref{Par.12}) with $C^\infty$ function $f=f(u,x)$; see
\cite{Pol02} and earlier references therein. The principal part of
the construction in \cite{Palis, Pol02} uses the fact that
the dimension  of manifold of equilibrium points
is less than the dimension of the corresponding kernel of the linearized operator
meaning a certain ``defect" of dimensions.

It is well-known that there are special (often similar to ours)
cases of branching of equilibria, where Hale--Raugel's approach
applies to guarantee convergence of the orbits of gradient
systems; see \cite{Hale92}, where a survey and  further references
are given. This approach essentially relies on spectral properties
of the linearized operator (main Hale--Raugel's hypothesis is that
0 has multiplicity at most one, or $k>1$ under special hypothesis)
and uses properties of stable, unstable and center manifolds, that
are difficult to justify for some less smooth equations. We also
refer to more recent paper \cite{Busca02}, where a similar
approach to stabilization is used and other references can be
found.

Our method is different and uses Zelenyak's ideas (1968) from one
dimension that are mainly connected with Lyapunov functions only.
It is important that Zelenyak's approach   yields the exponential
rate of convergence to a single equilibrium and can be extended to
non-autonomous perturbations of gradient systems, for which the
classic stable-centre manifold theory fails.



Let us  present the main ``resonance" hypothesis, under which
non-convergent orbits do not exist. We formulate it in a general
form, where, for the
 PDEs (\ref{Par.12}), by the Frechet
derivative we mean $V''(\psi)= \D + f'(\psi)I$ assumed to be a
Fredholm operator of index zero at any equilibrium point from $\Phi$.

\ssk

{\bf Hypothesis} ${\bf (H_R).}$ {\em If branching occurs at any
equilibrium $\psi$ and \, ${\rm dim \, ker} \, V''(\psi) =k \ge
1$, then the corresponding branching subset $S_k$ of equilibrium
points is a locally
smooth $k$-dimensional manifold $($a relatively open surface$)$}. 

\ssk

As a first comment, we note that, obviously, ${\bf (H_R)}$ always
holds in the {\em linear} case with $f(\psi)= \l \psi$, where
branching of equilibria  occurs for $ \l \in \s(\D)$. Then the
subsets coincide,
 \beq
 \label{LL1}
S_k={\rm ker} \, V''(\psi) \quad  (V''(\psi)= \D  + \l I).
 \eeq
Indeed in the linear case the convergence of all the orbits
follows from the completeness and closure of the eigenfunctions of
$\D$ as a straightforward consequence of general theory of
self-adjoint operators, \cite{BS}, and we do not treat this case
here. 

 The
condition in ${\bf (H_R)}$  cannot be improved in the sense that,
in a non-resonance branching, where the smooth  stationary
manifold $S_k$ has dimension $k < {\rm dim \, ker} \, V''(\psi)$,
this makes it possible to apply the idea of construction of
non-converging orbits from \cite{Palis, Pol02}.
According to our proof,
any  defect of dimensions  destroys the approach.
  Concerning the resonance branching
condition,
 we refer to the classical branching theory for
nonlinear equations in Banach spaces \cite[Ch. 7,8]{VainbergTr},
\cite[Ch.~6]{Berger}, \cite[\S~30.1]{Deiml} for the conditions,
which guarantee necessary branching in variational setting, and
also to \cite{Kiel}, where applications to elliptic problems of
interest are discussed in detail. Notice that, for (\ref{Par.11}),
 the branching parameter is not available explicitly in the
equation and such  branching cases are less studied in the literature.
 We will present simple examples showing that in standard  nonlinear elliptic problems,
the resonance branching with $k=N-1$ actually exists in dimensions
$N \ge 2$.

\subsection{Main result, plan and extensions}

 We now state the main result of the paper, which, as we have
 mentioned, are established by using key Zelenyak's ideas from
 \cite{Zel}.


 \begin{theorem}
 \label{Th.C1}
 Let $f$ be a $C^1$ function.
 Let $\Phi= \{\psi\}$ be the subset of all bounded  solutions
of the elliptic problem $(\ref{Par.11})$ and let ${\bf (H_{R})}$
hold. 
 Let, given bounded
initial data $u_0$, $u(x,t)$ be the unique classical solution of
$(\ref{Par.12})$. Then:

{\rm (i)} in the coercive case $(\ref{f.cond})$, there exists a
unique $\psi \in \Phi$ such that $\o(u_0) = \{\psi\}$;

{\rm (ii)} in the non-coercive case $(\ref{f.f.1})$, the same is
true for any uniformly bounded orbit;

{\rm (iii)} in both cases, the convergence to $\psi$ is exponential.

\end{theorem}

The hypothesis of boundedness of the orbit in (ii) is essential
since (\ref{Par.12})
 admits solutions, which blow up in finite time and hence exhibit entirely
 different asymptotic patterns (specific asymptotic techniques are explained
 in  \cite{AMGV}).
As we have mentioned,  for analytic nonlinearities $f(u)$ (for
(\ref{f.f.1}) this means
 $p=3,5,...$), the evolution completeness conclusion follows from
 the classical \L{ojasiewicz}--Simon inequality; see references in
 \cite{Chill, Fer}, where  versions of such
  inequalities were  proposed for some
non-analytic settings.  The result (iii) suggests that in the
resonance branching, the \L{ojasiewicz}--Simon inequality with
$\theta = \frac 12$ holds in a neighbourhood of $S_k$, though we
do not study this.

\ssk

 The layout of the paper is as follows. Theorem \ref{Th.C1} is
proved in Section \ref{SectSS}. In Section \ref{SectRes} we
present some examples showing  the possibility of resonance
branching. We
   show in Section \ref{SectRate}  that, for some $C^\infty$
 nonlinearities, the \L{ojasiewicz}--Simon inequality fails for any arbitrarily
 small $\theta >0$, so
 the stabilization  technique needs other non-rational moduli of continuity in (\ref{LS1}),
 which are proposed. It
 follows from (iii) in Theorem \ref{Th.C1} that such
 essentially non-analytic
 counterexamples correspond to isolated equilibria.

The results are true  for smooth higher-order
equations. E.g., Theorem \ref{Th.C1} is valid for a general
$2m$th-order parabolic flow for any $m \ge 1$
 \beq
 \label{2m.11}
 u_t = (-1)^{m+1} \D^m u + f(u) \quad \mbox{in} \,\,\, \O \times
 \re_+,
 \eeq
 with $m$ zero Dirichlet boundary conditions.
  The critical Sobolev exponent for  (\ref{f.f.1}) is
  \beq
  \label{pss}
  \mbox{$
  p_S= \frac{N+2m}{N-2m} \quad (N>2m). $}
  \eeq
We can also consider quasilinear fourth-order uniformly parabolic
equations like
 \beq
 \label{gg.1}
 u_t = - \D \, g(\D u) + f(u), \quad g'(s) \ge c_0>0,
  \eeq
or with others potential operators in various metrics. In
Section \ref{SectHi} we extend the approach to smooth gradient
systems in Hilbert spaces. In Section \ref{SectPert} we treat the
case of asymptotically small {\em non-autonomous} perturbations of
gradient systems. We show that the main result of convergence
(\ref{evC}) of bounded orbits remains true for perturbed parabolic
equations (\ref{Par.11}) such as
 \beq
 \label{eq.PP}
 u_t = \D u + f(u) + h(t) g(u), \quad \mbox{where} \,\,\, h(t) \to
 0 \,\,\, \mbox{as} \,\,\, t \to \infty,
  \eeq
  and $g(u)$ is a given smooth function, and also applies for
  similar perturbations of PDEs (\ref{2m.11}) and (\ref{gg.1}).
We then impose a certain restriction on the rate of decay of
$h(t)$ as $t \to \infty$, but, anyway, we can treat
non-exponentially small perturbations, e.g., $h(t) = O(\frac
1{t^3}$).

\section{\bf Proof of Theorem \ref{Th.C1}}
 \label{SectSS}

Thus we  consider the uniformly elliptic equation (\ref{Par.11})
and the corresponding uniformly parabolic flow (\ref{Par.12}). Notice that
the subset of all stationary solutions $\Phi=\{\psi\}$
includes  the trivial one $\psi=0$. Obviously, 0 is exponentially
stable for (\ref{Par.12}) with nonlinearity (\ref{f.f.1}) and
hence must be taken into account in the evolution completeness
analysis.

Let us mention again that according to the potential structure of the equation (\ref{Par.11})
and to (${\bf H_R}$), we are going to use the following properties:

\ssk

 \noi(i) $S_k$ is open, so any $\psi \in S_k$ has a
neighbourhood relatively open in $S_k$, and

\ssk

 \noi (ii) $V(\psi)=B={\rm const.},$ for all $\psi \in S_k$.

\subsection{Proof Theorem \ref{Th.C1}}

Thus we consider a smooth flow (\ref{Par.12}) in, say, $H_0^2(\O)$ and study the $\o$-limit set $\o(u_0)$ of
a uniformly bounded orbit $\{u(t), \, t>0\}$.
 Take an arbitrary equilibrium $\psi$ from the $\o$-limit set $\o(u_0)$, which is
non-empty, compact and connected for such smooth gradient systems,
\cite{Ha}. So that there exists a monotone sequence $\{t_k\} \to \infty$ such that
 \beq
\label{se1}
u(t_k) \to \psi \quad \mbox{as} \,\,\, k \to \infty.
 \eeq

First of all, if $0 \not \in \s(V''(\psi))$, where
 $V''(\psi)=\D + f'(\psi) I$, then
$\psi$ is isolated \cite{VainbergTr, KrasZ} and the convergence result (\ref{evC})
follows for given $u_0$. Second, assume that, under hypothesis (${\bf H_{R}}$), $0 \in
\s(V''(\psi))$ and $\psi$ is not isolated and belongs to a smooth
manifold $S_k$ branching at $\psi$, where $k$ coincides with the dimension of ${\rm
ker} \, V''(\psi)$.

Our further analysis uses the structure  of Zelenyak's
proof in the case $N=1$, \cite{Zel} and consists of seven steps.

\ssk

\noi{\bf Step 1: optimal approximation of $u(t)$ on $S_k$.} Let
$S_k = \{\psi_\mu\}$, where $\mu \in \re^k$ denotes  local
coordinates on $S_k$ such that $\psi=\psi_0$. Given a point $u(t)$
for some $t \gg 1$, when $u(t) \approx \psi$, we define the
optimal approximation $\psi_{\mu(t)}$ of $u(t)$ on $S_k$ by
minimizing the distance
 \beq
 \label{d12}
h(\mu)= \| u(t) - \psi_\mu\|,
 \eeq
 where $\| \cdot \|$ denotes the $L^2(\O)$-norm induced by the scalar
 product $\langle \cdot, \cdot \rangle$. Then ${\rm inf}_\mu
 h(\mu)$ is  attained at some $\mu=\mu(t) \approx 0$, and, by
 construction, the orthogonality condition holds,
  \beq
  \label{ort1}
  u(t) - \psi_{\mu(t)} \,\, \bot  \,\, {\rm ker} \, V''(\psi_{\mu(t)}).
   \eeq
In the further evolution analysis we
deal with equilibria $\psi_\mu$ in $S_k$ for $\mu \approx 0$.
 It suffices to prove that there exists the limit
 $$
\mu(t) \to 0 \quad \mbox{as} \,\,\, t \to \infty.
 $$


\ssk

  \noi{\bf Step 2: a priori bound for the linearized stationary equation.}
 We linearize the stationary equation (\ref{Par.11})
at $\psi=\psi_{\mu(t)}$ and consider the corresponding inhomogeneous problem
 \beq
\label{in1} V''(\psi_{\mu(t)}) w= g,
 \eeq
where $g \in L^2$ is a given function. Then, in view of the
orthogonality to the kernel, (\ref{ort1}), by the standard theory
of elliptic self-adjoint operators (see e.g. \cite{Berger, BS}),
there exists a constant $C_1>0$ such that
 \beq
\label{w12} \|w\|_{H^2} \equiv \|u(t) - \psi_{\mu(t)}\|_{H^2} \le
C_1 \|g\|.
 \eeq

\ssk

 \noi{\bf Step 3: spectral gap.} By hypothesis $({\bf H_{R}})$, $0 \in \s(V''(\psi_{\mu}))$
for all $\mu \approx 0$. Moreover, the resonance condition also assumes that the spectral gap
 \beq
\label{gg1} \Lambda(\mu)= \,\, \inf \{|\l|: \,\,\l \in \s({\bf
A}'(\psi_{\mu})), \, \l \not = 0\}
 \eeq
is uniformly bounded away from zero, i.e., there exists a constant $c_2>0$ such that
 \beq
\label{gg11}
\Lambda(\mu) \ge c_2>0 \quad \mbox{for all} \,\,\, \mu \approx 0.
 \eeq
This follows from the fact that, by $({\bf H_R})$, the kernel ${\rm ker} \, {\bf
A}'(\psi_\mu)$ changes continuously with $\mu \approx 0$ as the
tangent space of the changing continuously points $\psi_\mu$ on the given
smooth manifold $S_k$.


\ssk

 \noi{\bf Step 4: estimate for the linearized parabolic equation.}
 We now set  $w(t)=u(t)- \psi_{\mu(t)}$.
Consider next  the parabolic equation (\ref{Par.12}), linearize
the right-hand side at $u=w_{\mu(t)}$ and write down it as the
inhomogeneous elliptic equation (\ref{in1}). Then, by
linearization, we obtain an extra quadratic term $g_1(u)= \frac 12
\, f''(\xi_\mu)(u- \psi_{\mu(t)})^2$, so that, in (\ref{in1}),
 \beq
\label{F11} g= u_t + g_1(u), \quad \mbox{with} \,\,\, |g_1(u)| \le
C_3 \|u(t)- \psi_{\mu(t)}\|^2_{H^2}.
 \eeq
By convergence in such a smooth parabolic equation, we have that
the smallness of $w(x,t)=u(x,t)-w_{\mu(t)}(x)$ implies the
smallness of the derivatives, so that we may assume that
 $C_3$ is uniformly bounded on any intervals $[t_k,t_k+T]$
 as $k \to \infty$.

Combining the elliptic estimate (\ref{w12}) and the parabolic one
(\ref{F11})  yields the following bound:
 \beq
\label{w31}
 \|u(t) - \psi_{\mu(t)}\|^2_{H^2} \le C_4 \|u_t\|^2,
 \eeq
for all $t \gg 1$ such that the orthogonality (\ref{ort1}) holds.


\ssk

\noi{\bf Step 5: local exponential convergence of the Lyapunov function.}
The Lyapunov function $V(u(t))$ is strictly monotone increasing on evolution orbits,
 \beq
\label{V11} \mbox{$ \frac {{\mathrm d}}{{\mathrm d}t} \, V(u(t)) =
\int\limits_{\O} (u_t)^2 >0 \quad (u_t \not = 0), $} \eeq so that there
exists the finite limit
 \beq
\label{V12}
V(u(t)) \to B^+ \quad \mbox{as} \,\,\, t \to \infty,
 \eeq
and then $V \equiv B$ on $S_k$.
Fixing the parameter of the optimal approximation $\mu=\mu(t)$, by using
 the standard expansion at $\psi_\mu$, one can see
that
 \beq
\label{V13}
|V(u)- V(\psi_\mu)| \le C_5 \|u-\psi_\mu\|^2_{H^1} \le C_5 \|u-\psi_\mu\|^2_{H^2}.
 \eeq
Therefore from (\ref{V11}), by (\ref{w31}) and (\ref{V13}), we
have that
 \beq
\label{V15}
 \mbox{$ \frac {{\mathrm d}}{{\mathrm d}t} \, [B-
V(u(t))] = -\int\limits_{\O} (u_t)^2 \le - C_6 \|u- \psi_{\mu(t)}\|^2_{H^2} \le -
C_7[B- V(u(t))]. $} \eeq Integrating this inequality yields the
local exponential convergence of $V(u(t))$ to $B$ on any
arbitrarily large bounded intervals $t \in [t_k, t_k+T]$ as $k \to
\infty$,
 \beq
\label{V65}
B- V(u(t)) \le [B- V(u(t_k))] {\mathrm e}^{-C_7(t-t_k)}.
 \eeq

\ssk

\noi{\bf Step 6: local exponential convergence of $V(u(t))$
implies exponential convergence of $u(t)$ to equilibrium.} Here we
use Lemma 4 in \cite{Zel}, actually establishing a weighted
Gronwall's-type inequality by using a discrete partition
technique. Later on, finite partitions have been widely used for
more general   Gronwall's weighted inequalities; see Henry's
famous book \cite[Ch.~7]{He}.

\begin{lemma}
\label{Pr.Z4}
Let, for all $t >0$,
 \beq
\label{V16}
\mbox{$
\int\limits_t^\infty \int\limits_\O (u_t)^2 \le C_8 {\mathrm e}^{-t}.
 $}
 \eeq
Then there exists a constant $C_9>0$ 
  such
that
 \beq
\label{V17} \|u(t)-u(\t)\| \le C_9 {\mathrm e}^{-\frac t2} \quad
\mbox{for all} \,\,\, 0 < t \le \t < \infty.
 \eeq
\end{lemma}

 \noi {\em Proof.}  
Firstly, if $|t-\t| \le 1$ we apply the H\"older inequality to get
that
 $$
\begin{matrix}
\|u(t)-u(\t)\|^2
 = \int\limits_\O \big(\int\limits_t^\t u_t \big)^2 \le
 (\t-t)\int\limits_\O \big(\int\limits_t^\t (u_t)^2\big)
\le C_8 \, {\mathrm e}^{-t}.
 \end{matrix}
 $$
If $|t-\t| >1$, then we perform a uniform partition of the
interval $[t,\t]$ into $K$ parts of length 1 with a reminder, and
apply the H\"older inequality in each subinterval to get
 \beq
 \label{uu1}
\begin{matrix}
\|u(t)-u(\t)\|^2 \le \sum\limits_{i=0}^{K-1}
 \int\limits_\O \int\limits_{t+i}^{t+i+1} (u_t)^2
 + \int\limits_\O \int\limits_{t+i+K}^\t (u_t)^2
 \cr\cr
\le  C_8 \sum\limits_{i=0}^{K-1} {\mathrm e}^{-(t+i)} + C_8{\mathrm
e}^{-(t+i+K)} \le C_{9} {\mathrm e}^{-t}.
 \end{matrix}
 \eeq
 $\qed$

\ssk

\noi{\bf Step 7:  exponential estimate implies convergence to
$\psi$.} We continue to describe the evolution of $u(t)$ on a
large finite interval $[t_k,t_k+T]$ with $k \gg 1$, on which, by
assumption, the orbit enjoys all the estimates following from the
optimal approximation on the stationary subset $S_k$. So for
any  $ t \in [t_k,t_k+T]$, there holds:

\ssk
 (a) By (\ref{V65}) and (\ref{V11}),
 $$
\mbox{$ \int\limits_t ^\infty  \int\limits_\O (u_t)^2 = B-V(u(t))\quad \mbox{is exponentially small;}
 $}
 $$

 \ssk

 (b) By Lemma \ref{Pr.Z4}, 
  $\|u(t_k) - u(t)\|$ is exponentially small;

  \ssk

 (c) Therefore, $u(t) \approx u(t_k) \approx \psi$ for all $t \ge
 t_k$,
  completing the proof in the case (i). 

  \ssk

The case (ii) is the same once we have fixed a uniformly bounded
orbit.

\ssk

(iii) See (b) above. $\qed$

\section{\bf Discussion: resonance branching and exponential convergence}
\label{SectRes}

\subsection{Example: non-isolated equilibria and resonance branching}

  We take
the equation with the cubic analytic nonlinearity in the unit ball in $\ren$
 \beq
\label{D1}
\D \psi + \psi^3=0 \quad \mbox{in} \,\,\, B_1 = \{|x|<1\} \subset \ren, \quad \psi|_{\partial B_1}=0.
 \eeq
Then $p=3$ is in the subcritical Sobolev range if $3 <
\frac{N+2}{N-2}$, i.e., we need
 $
 N < 4.
 $
 Taking other $p<p_S$ provides us with similar examples for any $N
 \ge 2$.

Firstly, we consider  (\ref{D1}) in the half of the ball, $\O_+ =
B_1 \cap \{x_1 > 0\}$.  There exists a unique strictly positive
classical solution $\psi_+ \in H_0^1(\O_+)$ of (\ref{D1}) in
$\O_+$ constructed by the standard  variational technique,
\cite{Berger,PohFM}.

Next, since equation (\ref{D1}) is invariant under the
reflection $x_1 \mapsto -x_1$, setting $\psi= \psi_+$ in
$\O_+$ and $\psi=- \psi_+$ in $\O_-= B_1 \setminus \bar \O_+$ yields a non-radial
$H^1_0$-solution  of (\ref{D1}) in $B_1$.
 This extended function is a weak solution of (\ref{D1})
in the sense that it satisfies the corresponding integral
identity and hence is a classical solutions by the theory of uniformly elliptic equations.
 Any smooth
invariant orthogonal transformation in $\ren$ 
leaving the Laplacian and the Dirichlet boundary condition
invariant produces continuous families of solutions consisting of
non-isolated points.  By the classical branching
theory \cite{KrasZ},  
 in this
case, $0$ belongs the spectrum of the corresponding  linearized
operator $V''(\psi)= \D + 3 \psi^2 I$.

 In particular, we fix $N-1$ angels of
rotations $\{\theta_1,...,\theta_{N-1}\}$ of the polar coordinate
system in $\ren$ to get at least  $N-1$ linearly independent
elements of the kernel of 
$\D + 3 \psi^2 I$ given by
 $$
 \mbox{$
 \varphi_k= \frac{{\mathrm d}}{{\mathrm d} \theta_k} \, \psi, \quad
 k=1,...,N-1.
 $}
  $$
  Bearing in mind these rotations,
 the
stationary subset $S_{N-1}$ is expected to be precisely ({\em
N}$-$1)-dimensional, and this  $k=N-1$ coincides with the
kernel dimension. Unfortunately, we cannot prove the exact
equality, and present later on another non-local model, for which
Hypothesis (${\bf H_R}$) is guaranteed.
 For the problem on the plane ($N=2$), the
corresponding eigenfunction of $V''(\psi)$ with $\l =0$
 is given by rotation by the single angle $\theta=\theta_1$, so
  \beq
\label{psi1NN}
\mbox{$
 \varphi_1 = \frac{{\mathrm d}}{{\mathrm d} \theta} \psi,
$}
 \eeq
 and  the kernel seems to be one-dimensional (we do not prove
  this either).

Using similar reflections, one can construct in the unit ball $B_1
\subset \ren$ the nonlinear eigenfunctions $\psi_l(x)$ satisfying
(\ref{D1}) with arbitrarily number $l \ge 1$ of connected
positivity and negativity components. By orthogonal invariant
transformation, each one generates continuous families of other
stationary solutions.

 \subsection{On exponential convergence}

We first illustrate the reason for the exponential convergence by using the above example
in $\re^2$, where the kernel is one-dimensional with the basis function (\ref{psi1NN}).


\ssk

 \noi {\bf Orthogonality conditions for  branching.} The branching occurs in (\ref{D1}) from
solution $\psi_0=\psi$, so by the classical branching theory \cite{VainbergTr}, we look for
a smooth curve of solutions in the form
 \beq
\label{psimu}
\psi_\mu = \psi_0 + \mu \psi_1 + \mu^2 \psi_2 + \mu^3 \psi_3+ \mu^4 \psi_4+... \, ,
 \eeq
where, up to scaling, $\mu$ can be attributed to the angle of
rotation. Substituting (\ref{psimu}) into (\ref{D1}) and equating
similar terms yields the system for expansion coefficients
 \beq
\label{sys1}
\left\{ \begin{matrix}
 \D \psi_0 + \psi_0^3=0,\, \qquad \qquad \qquad \quad \,\qquad \qquad \qquad  \cr
 \D \psi_1 + 3 \psi_0^2 \psi_1=0, \,\qquad \qquad \qquad \qquad \qquad \quad\,\,\, \cr
\D \psi_2 + 3 \psi_0^2 \psi_2=- 3 \psi_0 \psi_1^2,\, \qquad \,\,\,
\qquad \qquad \qquad\cr
  \D \psi_3 + 3 \psi_0^2 \psi_3=- \psi_1^3 -
6 \psi_0 \psi_1 \psi_2,  \qquad \qquad \quad \,\, \cr
 \D \psi_4 +
3 \psi_0^2 \psi_4=- 3 \psi_0 \psi_2^2 - 3  \psi_1^2 \psi_2 - 6
\psi_0\psi_1 \psi_3, \cr ... \qquad ... \quad \quad ... \,\, .
\end{matrix}
 \right.
 \eeq
The second equations  says that $\psi_1$ is from the non-empty
kernel and hence is given by (\ref{psi1NN}). The third equation
then yields  the first orthogonality condition
 \beq
\label{or77}
\mbox{$
 \int\limits_\O \psi_0 \psi_1^3 = 0,
$}
 \eeq
which is necessary for branching to occur.
Determining  $\psi_2$ up to a constant $\b_2$,
 $$
\psi_2 = \bar \psi_2 + \b_2\psi_1, $$ where $\bar \psi_2$ is a
solutions of the inhomogeneous equation, and substituting into the
fourth one, for its solvability we have to have that  
 $$
\mbox{$
 \int\limits_\O \psi_1^4 +  6 \int\limits_\O \psi_0 \psi_1^2 \bar \psi_2 + 6 \b_2 \int\limits_\O
\psi_0 \psi_1^3
=0.
$}
 $$
Since the last term containing $\b_2$ vanishes by (\ref{or77}),  we arrive at the
second orthogonality condition
 \beq
\label{or22} \mbox{$ \int\limits_\O \psi_1^4 +  6 \int\limits_\O \psi_0 \psi_1^2
\bar \psi_2=0.
 $}
 \eeq
Next choosing similarly
 $$
\psi_3 = \bar \psi_3 + \b_3 \psi_1,
$$
we obtain from the fifth equation that
 $$
\mbox{$ 2 \int\limits_\O \psi_0 \psi_1^2\bar \psi_3 +  \int\limits_\O \psi_0
\psi_1 \bar \psi_2^2 + \int\limits_\O \psi_1^3 \bar \psi_2 + \frac 23 \,
\b_2 \int\limits_\O \psi_1^4=0. $}
 $$
This determines $\b_2$
and so on.
Thus, in general, this kind of branching demands more than one
orthogonality condition (and, of course,
 existence of a non-trivial kernel).




\ssk

 \noi {\bf On a centre manifold link.}
Let us now discuss how these branching  conditions affect the
evolution properties
 of the corresponding parabolic equation
 \beq
\label{Pa33}
u_t = \D u + u^3 \quad \mbox{in} \,\,\, B_1 \times \re_+.
 \eeq
In this analytic case, 
there exists a local centre $C^\infty$ manifold $W_{\rm loc}^{\rm
c}(\psi_0)$; see general theory in \cite{Pazy}, \cite[\S~9]{Lun},
\cite{Sim1} for applications to reaction-diffusion equations and
\cite{Hale92}, where the convergence problem is studied.

 Since  $W_{\rm loc}^{\rm
c}(\psi_0)$ is tangent to the corresponding centre subspace
$E^{\rm c}(0)= {\rm Span}\{ \psi_1\}$, looking for the centre
manifold behaviour with solutions of the form (see references and
accurate estimates for the non-analytic flow in the next section)
 $$
u(x,t) = \psi_0(x) + a_1(t) \phi_1(x)+... \quad \mbox{for} \,\,\, t \gg 1,
 $$
and projecting the PDE (\ref{Pa33}) on $E^c(0)$, we obtain a perturbed ODE
 \beq
\label{Od11}
a_1'= \g a_1^2 + ... \, ,
 \eeq
where, as one can see, the constant $\g$ vanishes,
 \beq
\label{gg1uu}
\mbox{$
 \g = \int\limits_\O \psi_0 \psi_1^3 =0,
$}
 \eeq
by the first orthogonality branching condition (\ref{or77}).
Therefore the higher-order terms should be taken into account, but
 these seem also do not help to detect a suitable non-trivial evolution on the centre manifold.

Indeed one can see that this one-dimensional branching manifold
$S_1$ is a centre manifold for the parabolic problem. Under the
above assumptions, denoting by $\{\phi_k\}$ the complete subset of
eigenfunctions of $V''(\psi_0)= \D + 3 \psi_0^2 I$, we may assume
the following expansion (cf. (\ref{psimu}))
 $$
S_1 = \{ \psi_\mu = \psi_0 + a_1(\mu) \phi_1 + a_2(\mu) \phi_2+...\},
 $$
where $a_k(0)=0$ for any $k$. Therefore, the evolution on this
centre manifold with 
 $$
\mbox{$ u(x,t) =\psi_0(x) + \sum\limits_{(k)} a_k(t) \phi_k(x) $}
 $$
 is governed by the trivial dynamical system
 $$
 a_k'=0, \quad k=1,2,... \, ,
 $$
which we have observed in (\ref{Od11}), (\ref{gg1uu}). Notice that
a centre manifold consisting of pure equilibria is not an
exceptional situation in parabolic problems; cf. \cite{Sim2},
where such an invariant exponentially stable centre manifold has
been detected in the free-boundary Mullins--Sekerka model (a kind
of Hele--Shaw flow).

Thus in this case, as well as and in other resonance cases, there
exists {a local centre manifold for the parabolic problem
precisely coinciding with the stationary subset $S_k$}, so no
evolution of orbits on $W^{\rm c}_{\rm loc}(\psi_0)$ can be
observed. Then the exponential convergence  in the Theorem
\ref{Th.C1} (iii) can be associated with  the fact that, in the
orthogonal complement of the kernel, the behaviour is purely
exponential and corresponds to the evolution on the orthogonal
stable manifold constructed at a different equilibrium
$\psi=\psi_{\mu(t)} \approx \psi_0$; cf. Step 1 of the proof in
Section \ref{SectSS}. In the invariant manifold theory, this
 is usually expressed by the fact that, for
$C^2$ nonlinearities and good linearized sectorial operators with
discrete spectrum and finite-dimensional unstable and centre
subspaces, etc., the centre manifold is exponentially stable; see
e.g. \cite[Prop.~9.2.3]{Lun} and \cite{Sim1}.


\section{\bf An explicit example with non-local nonlinearity}

Let  $\O$ be a bounded smooth domain in $\ren$
  and  let $\{\l_k\}$ and $\{\psi_k\}$
be the eigenvalues of $\D$, where
 each one $\l_k$ repeated as many
times as its multiplicity $\kappa_k \ge 1$,
 and the corresponding complete, closed subset of  eigenfunctions that are
 orthonormal in $L^2(\O)$.
The case of the unit ball, $\O=B_1$, is classical. Here, for $\D$,
all the eigenvalues, eigenfunctions, multiplicities, etc., are
described by
 the Laplace-Beltrami operator $\D_\s$ on the unit sphere $S^{N-1}=
\partial B_1$ in $\ren$, which
 in the polar coordinates $(r,\s)$ is given by
\beq
 \label{L-B}
  \D= \D_r + \mbox{$\frac 1{r^2}$} \,\D_\s,
 \eeq
$\D_\s$
  is a regular operator
with the discrete spectrum in $L^2(S^{N-1})$ (again each one
repeated as many times as its multiplicity)
 \beq \label{LBsp} \s(-\D_\s) =
\{\nu_k=k(k+N-2), \,\, k \ge 0\}
 \eeq
and an orthonormal, complete,  closed subset $\{V_k(\s)\}$ of
eigenfunctions, which are homogeneous harmonic $k$-th order
polynomials restricted to $S^{N-1}$.

Fix an $l \ge 2$  and 
  consider the PDE with a
non-local cubic nonlinear term
 \beq
 \label{cu1}
 \mbox{$
 u_t = \D u + \l_l u - \big(\int\limits_\O u^2\big) u
 $}
 \eeq
 with the same Dirichlet boundary conditions and bounded initial
 data $u_0$. Such non-local
 parabolic models can provide us with further examples of delicate
 asymptotics obtained via explicit computations; see \cite{PohFM}.
Firstly, it is easy to solve the stationary  equation
 \beq
 \label{ellp1}
  \mbox{$
 \D \psi + \l_l \psi - \big(\int\limits_\O \psi^2 \big) \psi=0.
  $}
  \eeq
  Namely, there exist non-trivial equilibria of the form
   \beq
   \label{ps1}
   \psi(x) = \pm \sqrt{\l_l-\l_j} \, \psi_j(x) \quad \mbox{provided
   that} \,\,\, \l_l > \l_j,
    \eeq
    and other obviously  constructed linear combinations of such functions
    corresponding to the same eigenvalues.

 \subsection{Result on convergence}

Secondly, studying the parabolic equation (\ref{cu1}) and
  looking for solutions in the form
 of the eigenfunction expansion
  $$
  \mbox{$
  u(x,t)= \sum\limits_{(k)} a_k(t) \psi_k(x),
   $}
   $$
we obtain the following dynamical system for the expansion
coefficients:
 $$
 \mbox{$
  a_k'=(\l_l-\l_k- |a|^2)a_k, \,\,\, k=1,2,... \, , \quad
  \mbox{where} \,\,\,
  |a|^2= \sum a_k^2.
   $}
  $$
Integrating yields
 $
 a_k(t) = a_k(0){\mathrm e}^{(\l_l-\l_k)t}{\mathrm e}^{-\int\limits_0^t |a|^2(s) \, {\mathrm d}
 s}.
 $
 Calculating the sum $|a|^2$,
  we  derive 
   $$
   \mbox{$
   |a|^2= \Psi(t){\mathrm e}^{-\int\limits_0^t |a|^2(s) \, {\mathrm d}
 s}, \quad \mbox{with} \,\,\, \Psi(t)= \sum {\mathrm
 e}^{2(\l_l-\l_k)t}a_k^2(0).
 $}
  $$
Setting $Z={\mathrm e}^{-\int_0^t |a|^2(s) \, {\mathrm d}
 s}$ yields a simple ODE,
  $
  Z'= - 2 \Psi(t) Z^2,
  $
  so, on integration, 
 $$
 \mbox{$
 |a|^2(t)=\frac {a_l^2(0)+\sum\limits_{\l_m \not = \l_l}{\mathrm
 e}^{(\l_l-\l_m)t}a_m^2(0)
 }
 {1+ 2 a_l^2(0)t +
 \sum\limits_{\l_m \not =\l_ l} \frac 1{\l_l-\l_m}({\mathrm
 e}^{2(\l_l-\l_m)t}-1)a_m^2(0)  } \, ,
  $}
  $$
  where for simplicity  $a^2_l(0)$ denotes the sum
   $
   \mbox{$
    \sum\limits_{\l_k=\l_l} a^2_k(0).
   $}
   $
   Finally we arrive at the following expressions for the
   expansion coefficients:
 \beq
 \label{ak22}
 \mbox{$
 a_k(t)=a_k(0)\, \frac { {\mathrm e}^{(\l_l-\l_k)t}}{\sqrt{1+  2
 a_l^2(0)t+
 \sum\limits_{\l_m \not = \l_l} \frac 1{\l_l-\l_m}({\mathrm
 e}^{2(\l_l-\l_m)t}-1)a_m^2(0)}} \, .
  $}
  \eeq
This
 reveals three different
cases of solutions with  exponential and algebraic rate of
convergence to equilibria, depending on the dominant term in the
long square root in the denominator.
Namely, we take initial function
 $$
 \mbox{$
 u_0(x) = \sum a_k(0) \psi_k(x)
 $}
 $$
 such that
 there exists $j > 1$, for which
 \beq
 \label{k01}
 a_1(0)=...=a_{j-1}(0)=0 \quad \mbox{and} \quad a_j(0) \not = 0.
  \eeq
It is easy to derive from (\ref{ak22}) the following result.

\begin{proposition}
\label{Pr.NL}
 Let $(\ref{k01})$ hold. Then, as $t \to \infty$,
  \beq
  \label{u11}
{\rm (i)} \quad  u(x,t) = O({\mathrm e}^{-(\l_j-\l_l)t}) \to 0, \quad
\mbox{if \,$\l_j >
  \l_l$};
 \eeq
  \beq
  \label{u111}
  \mbox{$
{\rm (ii)} \quad   u(x,t)= O\big(\frac 1{\sqrt t}\big) \to 0, \quad \mbox{if \,
$\l_j =
  \l_l$}; \quad \mbox{and} \quad
   $}
   \eeq
  \beq
  \label{u1111}
{\rm (iii)} \quad  u(x,t)= \psi(x) + O(t \, {\mathrm
e}^{-2(\l_l-\l_j)t}),
   \quad
\mbox{if \, $\l_j <
  \l_l$},
   \eeq
   where $\psi \not =0$ is an equilibrium.
   \end{proposition}

It is important that the case (iii), where branching of equilibria
is available by (\ref{ps1}), precisely demonstrates that the
convergence to non-trivial stationary solutions is  exponential;
cf. Theorem \ref{Th.C1}(iii). It is worth mentioning
that 
the rate of convergence in (\ref{u1111}) is not purely exponential
and, in general, contains a lower-order algebraic factor $t$ that
is induced by the $O(t)$-term in the square root in (\ref{ak22}).
Different $\ln t$-perturbations can occur in problems with usual
(local) nonlinearities; see the Remark in Section \ref{SectRate}.

 The only case in Proposition \ref{Pr.NL}, where the convergence
 is not exponentially fast, is (ii), in which, by (\ref{ps1}),
the only equilibrium is trivial, $\psi=0$, and is isolated.

\subsection{Hypothesis (${\bf H_R}$) is valid}

We now present a rigorous evidence that Hypothesis $({\bf H_R})$
is valid for such operators. Consider the elliptic equation
(\ref{ellp1}). Fix a $j>l$, let $\kappa_j =1+k \ge 2$ be the
multiplicity of $\l_j$ and let $\{\psi_{j,1},...,\psi_{j,k+1}\}$
be the orthonormal subset of eigenfunctions of $\D$ corresponding
to $\l_j$. We fix an equilibrium
 \beq
 \label{eq.1}
 \mbox{$
 \hat \psi(x) = \sum\limits_{i=1}^{k+1} \hat c_i \psi_{j,i}(x),
  $}
  \eeq
  where, on substitution into (\ref{ellp1}), the coefficients
  satisfy (cf. (\ref{ps1}))
  \beq
  \label{eq.2}
  \mbox{$
  \sum \hat c_i^2= \l_l- \l_j.
  $}
 \eeq

 \noi{\sc Stationary subset $S_k$.} Obviously the subset $S_k$
 containing given equilibrium (\ref{eq.1}) is described by
  \beq
  \label{eq.3}
  \psi=  \mbox{$
 \hat \psi(x) = \sum\limits_{i=1}^{k+1}  c_i \psi_{j,i}, \quad \mbox{where}
 \,\,\, \sum c_i^2= \l_l- \l_j.
  $}
  \eeq
Therefore
 \beq
  \label{eq.4}
  {\rm dim} \, S_k=k=\kappa_j-1 \ge 1.
   \eeq

  \noi {\sc Kernel of $V''(\hat \psi)$.} It follows that
  \beq
  \label{eq.5}
  \mbox{$
  V''(\hat \psi)v= \D v+ \l_l v - 2 \hat \psi\big(\int \hat \psi v\big) -
  v\big(\int \hat \psi^2\big).
   $}
   \eeq
Substituting (\ref{eq.1}) yields the following equation for the
kernel:
 \beq
  \label{eq.6}
  \mbox{$
  V''(\hat \psi)v \equiv  \D v+ \l_j v - 2 \hat \psi(\int \hat \psi v)=0.
   $}
   \eeq
Finally taking $v$ in the form
 $$
 \mbox{$
 v= \sum\limits_{i=1}^{k+1} b_i \psi_{j,i}
  $}
 $$
 and substituting into (\ref{eq.6}) gives a single condition on
 $k$+1 expansion coefficients
 $$
  \mbox{$
 \sum\limits _{i=1}^{k+1} b_i \hat c_i=0,
  $}
  $$
  so that the kernel is precisely $k$-dimensional.

  By (\ref{eq.4}) this completes
  the analysis and proves that $({\bf H_R})$ is always valid for
  such non-local elliptic operators. Without any changes
 this example extends to such operators of
   arbitrary order
   $$
   \mbox{$
   V'(\psi)= -(-\D)^m \psi + \l_l \psi-\big(\int\limits_\O \psi^2\big) \psi
   \quad \mbox{with any $m \ge 1$}.
   $}
   $$

\section{\bf Rate of convergence can be  arbitrarily slow: $\theta=0$}
\label{SectRate}

Here we present some estimates showing that the rate of
convergence in parabolic problems with specially designed
nonlinearities can be arbitrarily slow.


\subsection {The original semilinear parabolic equation}

Let $\O$ be a bounded smooth domain in $\ren$ such that $\l_1=-1$
is the first simple eigenvalue of $\D$ in $L^2(\O)$ with domain
$H^2_0(\O)$ and the normalized eigenfunction $\phi_1>0$ in $\O$.
Consider the following parabolic equation:
 \beq
 \label{eq2}
 u_t = V'(u) \equiv  \D u + u - {\mathrm e}^{- 1/ {u^2}} \quad \mbox{in} \,\,\,
  \O \times \re_+, \quad u=0 \quad \mbox{on} \,\,\, \partial \O \times
  \re_+,
  \eeq
where we set $f(0)=0$ by continuity.
  Here
 $ {\mathrm e}^{- 1/ {u^2}}$ is the standard $C^\infty$  function, which is not analytic at
  $u=0$. In view of  the perfect spectral properties of $\D$, we have
  that there exists a local invariant, $C^\infty$, one-dimensional centre manifold
  $W^{\rm c}_{\rm loc}(0)$ of operator
  $V'(u)$,
  which is tangent to the eigenspace $E^{\rm c}(0)={\rm Span}\{\phi_1\}$
of $\D + I$; see \cite[Thm.~9.2.2]{Lun} and \cite{Sim1}. Looking
for the corresponding center manifold behaviour, we   decompose
the solution in the form
 \beq
 \label{vv1}
 u(x,t) = a_1(t) \phi_1(x) + u_2(x,t) \quad \mbox{for} \,\,\, t \gg 1,
  \eeq
  where $u_2(x,t) = \sum\limits_{k \ge 2} a_k(t) \phi_k(x) \perp E^{\rm c}(0)$, $u_2(\cdot,t) = o(a_1(t))$ as $t \to
 \infty$
 and we assume that $a_1(t) > 0$ for  $t \gg 1$, e.g., we take positive solutions.
  Projecting the PDE onto $E^{\rm c}(0)$ yields
   \beq
   \label{Pr.1}
   a_1' = - \langle  {\mathrm e}^{- 1/ {u^2}}, \phi_1
   \rangle \equiv  - \langle  {\mathrm e}^{- 1/(a_1^2 \phi_1^2+...)}, \phi_1
   \rangle
   \eeq
   for $t \gg 1$,
where $\langle \cdot, \cdot \rangle$ denoted the inner product in
$L^2(\O)$. Denoting $\rho_1= \max \phi_1(x)>0$ and $\rho_2 = \int
\phi_1>0$ yields the ordinary differential  inequality
   \beq
   \label{In1}
   a_1' \ge - \rho_2 {\mathrm e}^{- 1/{4\rho_1^2}{a_1^2}},
   \eeq
   so that $a_1(t) \ge \bar a_1(t)$, where $\bar a_1$ solves the
   corresponding ODE with the equality sign in (\ref{In1}).
   Finally, we obtain the following estimate:
    \beq
    \label{a33}
     \mbox{$
    a_1(t) \ge \frac {1}{3\rho_1 \sqrt{ \ln t}} \quad \mbox{for} \,\,\, t
    \gg 1
     $}
     \eeq
 for the stabilization on the centre manifold.
Obviously, this corresponds to $\theta=0$ in (\ref{th1}), so that
the rational algebraic modulus $|\cdot|^{1-\theta}$ in (\ref{LS1})
cannot be applied to such $C^\infty$ nonlinearities. One can
introduce an appropriate modulus, for which the rate (\ref{a33})
is acceptable. For instance, for a slightly modified function
 $$
 \mbox{$
 V(u) = \frac 12 \, {\mathrm e}^{-1/u^2} \,\,\Longrightarrow \,\,
 V'(u)= \frac 1{u^3} \, {\mathrm e}^{-1/u^2},
  $}
 $$
we have the  {\em generalized  
gradient inequality} inequality at $\psi=0$ for $u \approx 0$
(actually, it is equality)
 $$
 \tilde \o(V(u)) \le |V'(u)| \quad \mbox{with modulus} \,\,\, \tilde \o(s) =
 2|s| \, |\ln(2s)|^{3/2},
 $$
 which is ``almost" linear, strictly concave function as $s \to 0$.

It is easy to present other $C^\infty$, non-analytic
nonlinearities generating arbitrarily slow rate of convergence.
Recall that, according to the results in Section \ref{SectRes},
any such slow convergence  corresponds to isolated equilibria.

\ssk

\noi{\bf Remark: a cubic nonlinearity.} In the presence of
non-trivial kernels, a full
asymptotic expansion of solutions can be 
 a difficult problem even in the analytic cubic case
 \beq
 \label{be1}
 u_t = \D u + u - u^3.
  \eeq
It was shown in \cite{Kond1} that, under the same kernel
assumption
  \beq
  \label{la1}
 {\rm dim \, ker} \,(\D + I)=1,
  \eeq
  there exist solutions with the following logarithmically
  perturbed decay as $t \to \infty$:
   $$
   \mbox{$
   u(x,t) = \sum\limits_{k=0}^n t^{-\frac 12-k} \sum\limits_{j=0}^k \varphi_{kj}(x)(\ln t)^j +
   O(t^{-\frac 32 -n}),
   $}
    $$
    for some integer $n>0$, where $\varphi_{kj}$ are solutions of
    well-posed elliptic problems. This corresponds to a
    special
    case of centre manifold behaviour, where the leading term of
    convergence to $\psi=0$ is of order
    $O(t^{-\frac 12})$ that does
    not contain the $\ln t$-factor. It is easy to see why  in
    this case branching of equilibria from zero is not possible. Indeed, if this
    occurs, then in view of (\ref{la1}), by branching theory \cite{VainbergTr},
      the equilibrium branch should have
    the representation $\psi = \mu \varphi_0+...$, where $\mu$ is
    the branching  parameter and $\varphi_0$ is the eigenfunction
    of $\D$ with $\l=-1$. Substituting this into the stationary
    equation $\D \psi + \psi - \psi^3=0$ and multiplying by
    $\varphi_0$ yields $\int_\O \varphi_0^4=0$, a contradiction,
    so $\psi=0$ is the isolated equilibrium.




\subsection{A non-local semilinear parabolic equation}

As usual the computations are simplified for semilinear non-local
parabolic flows like
 \beq
\label{NN.1}
\mbox{$
u_t = \D u + \mu u - g'\big(\int\limits_\O u^2\big) u,
$}
 \eeq
where $g'(s)$ is a given $C^\infty$, non-analytic function.
The potential here is
 \beq
\label{NN.2}
\mbox{$
V(u(t))= - \frac 12 \, \int |D u(t)|^2  + \frac \mu 2 \, \int u^2 - g\big(\int u^2\big).
$}
 \eeq
Assuming again that $\mu = \l_1$ is the first eigenvalue of $\D$
in $L^2(\O)$, we obtain that the centre manifold behaviour
(\ref{vv1}) is described by the ODE
 \beq
\label{NN.3}
a_1' = - g'(a_1^2+...) \quad \mbox{for} \,\,\, t \gg 1.
 \eeq
Choosing non-analytic $C^\infty$ functions $g$, e.g.,
 $$
g(s) = s^{3/2} {\mathrm e}^{-1/s} \quad \mbox{for} \,\,\, s>0,
 $$
integrating (\ref{NN.3}) asymptotically yields a non-algebraic
decay
 $$
 \mbox{$
 a_1(t) \approx  \frac 1{\sqrt {\ln 2t}} \quad \mbox{as} \,\,\, t \to \infty.
  $}
 $$
Notice that,
besides slow decay behaviour,
 the present ``less nonlinear" model (\ref{NN.1}) is suitable for revealing
refined evolution  properties of orbits describing stabilization phenomena.

\section{\bf On applications to smooth gradient systems in Hilbert spaces}
\label{SectHi}

Such an extension is straightforward, since, in the proof of Theorem
\ref{Th.C1} in Section \ref{SectSS}, we have minimally used  the
specific properties of the second-order  elliptic and parabolic
 equations under consideration.

 Denoting by  $\langle \cdot, \cdot \rangle$ and $\|\cdot \|$ the
inner product and the induced norm in a separable  Hilbert space $H$, we
consider a smooth gradient flow in $H$,
 \beq
\label{H1} u_t = V'(u) \quad \mbox{for} \,\,\, t>0,
 \eeq
where the operator $V'$ with a dense, compactly embedded domain
$H^2=D(V') \subset H$
 is the Frechet derivative of a
$C^2$-functional $V:H^2 \to \re$. We  assume that 
 $V''(u):H^2 \to H$ is a Fredholm operator of index zero
 admitting a suitable self-adjoint extension with discrete
 spectrum, compact resolvent and a subset of
 eigenfunctions, which is  complete and closed in $H$.
 We impose the
necessary condition of coercivity of $V'$ to guarantee existence
of global orbits $\{u(t)\}$ for arbitrary initial data $u_0 \in
H^2$, which are sufficiently regular to satisfy the gradient
identity (cf. (\ref{V11}))
 \beq
\label{V11N}
\mbox{$
\frac {{\mathrm d}}{{\mathrm d}t} \, V(u(t)) = \|u_t\|^2 \ge 0. 
$}
\eeq

We impose the same main Hypothesis (${\bf H_{R}}$) and will next
check, using the scheme of the proof  from Section \ref{SectSS},
which conditions we need to guarantee the result (i) and, hence,
(iii) in Theorem \ref{Th.C1}.

\ssk

\noi {\bf Steps 1--3.} The construction remains the  same and
(\ref{d12}), (\ref{w12}) and (\ref{gg11}) are guaranteed by
assumed good spectral properties of the linearized self-adjoint
operator $V''(\psi_\mu)$ for any $\psi_\mu \in S_k$ uniformly in
$\mu \approx 0$.

\ssk

\noi{\bf Step 4.} Estimate (\ref{w31}) is  valid for sufficiently
smooth operator $V'$ in $H^2$.

\ssk

 \noi{\bf Step 5.} Here we need  estimate (\ref{V13}),
which is indeed Lagrange's formula for the smooth functional $V$
in $H^2$.

\ssk

\noi{\bf Step 6.} Lemma \ref{Pr.Z4} remains valid in the topology
of $H$.

\ssk

 \noi{\bf Step 7} remains unchanged.

\ssk

\noi Hence  Theorem \ref{Th.C1} (i), (iii) is true for
arbitrary smooth gradient flows under the presence of resonance
branching. It follows that the convergence (evolution
completeness) result holds for the $2m$th-order PDEs (\ref{2m.11})
and for other smooth gradient parabolic  flows governed by higher-order
operators.

\section{\bf Non-autonomous perturbations of gradient systems}
\label{SectPert}

More modifications of the approach are necessary to prove the
convergence result in Theorem \ref{Th.C1} (i) for non-autonomous perturbations of
(\ref{H1}). The main features of our analysis are illustrated
by the following example:
 \beq
 \label{H1P}
 u_t = V'(u) + h(t) W'(u) \quad \mbox{for} \,\,\, t>0,
  \eeq
where $W(u):H^2 \to \re$ is a $C^1$-functional. Without loss of
generality and for simplification of future manipulations, we
assume that the decay rate of perturbation $h \in C^2$ satisfies
 \beq
 \label{h1}
 h(t) \to 0 \,\,\,\mbox{as} \,\,\, t \to \infty,
  \quad
 \mbox{and} \quad
 h(t)>0, \,\, h'(t) < 0 \quad \mbox{for} \,\,\, t \ge 0,
  \eeq
so that, in particular, $ h' \in L^1(\re_+)$.
Otherwise, for non-monotone and changing sign rates of
perturbations, we can use estimates from above and below
with functions $h_\pm(t)$ satisfying necessary assumptions.

Actually, we will need a more restrictive condition on the decay
rate, 
 \beq
\label{h12N}
\sqrt h \in L^1(\re_+).
   \eeq
It is important to deal with  non-exponentially small
perturbations such as
 \beq
 \label{h11}
 \mbox{$
 h(t) = \frac 1{(1+t)^{\a}}
\quad \mbox{with any} \,\,\, \a >2,
  $}
  \eeq
for which (\ref{h12N}) holds.
  For the future purpose of integration of ordinary differential
  inequalities, we characterize this class of
 such slow decaying functions as follows: for any
constant $\b > 0$,
 \beq
 \label{h111}
 \mbox{$
 \int\limits_0^t h(s){\mathrm e}^{\b s} \, {\mathrm d}s = \frac 1 \b \,
 h(t){\mathrm e}^{\b t}(1+o(1)) \quad \mbox{as} \,\,\, t \to
 \infty,
  $}
 \eeq
 or, equivalently, by L'Hospital's rule,
  \beq
  \label{h12}
  \mbox{$
  \frac{h'(t)}{h(t)} \to 0^- \quad \mbox{as} \,\,\, t \to \infty.
  $}
   \eeq

Since (\ref{H1P}) is not a gradient system in the sense that a
monotone Lyapunov function, in general, does not exist, we need
first to prove that the crucial characterization (\ref{Om1}) of the $\o$-limit set
remains valid, so we should begin with

\ssk

\noi{\bf Step 0: $\o(u_0) \subseteq \Phi$.} Multiplying
(\ref{H1P}) by $u_t$ in $H$, instead of (\ref{V11N}), we obtain
the identity
 \beq
 \label{V11NP}
 \mbox{$
\frac {{\mathrm d}}{{\mathrm d}t} \, V(u(t)) + \frac {{\mathrm
d}}{{\mathrm d}t} \,(h(t) W(u)) - h'(t) W(u)
 = \|u_t\|^2 \ge 0. 
$} \eeq
 Integrating over $(t,\infty)$ and using that $h' \in L^1(\re_+)$ yields the convergence
 \beq
 \label{ut1}
 \mbox{$
 \int\limits_t^\infty \|u_t(s)\|^2 \, {\mathrm d}s < \infty.
 $}
  \eeq
  Hence,
given a sequence $\{t_k\} \to \infty$ such that $u(t_k) \to \psi
\in \o(u_0)$, we obtain by the H\"older inequality that, for
 arbitrarily large fixed $t>0$,
 \beq
 \label{ut2}
 \mbox{$
 \|u(t_k+t)-u(t_k)\|^2 \le t \,  \int\limits_{t_k}^\infty \|u_t\|^2 \, {\mathrm d}s
\to 0 \quad \mbox{as}  \,\,\, k \to \infty.
 $}
 \eeq
 Passing to the limit as $t=t_k+t \to \infty$ in (\ref{H1P}), we then conclude
 that $u(t_k+t)$ converges in $H$ to a solution $u(t)$ of the
 limit autonomous equation (\ref{H1}), which is independent of
 $t$, so it is an equilibrium  $\psi \in \Phi$.

 We now return to seven steps of the proof in Section
 \ref{SectSS}.

\ssk

 \noi{\bf Steps 1-3} are unchanged.

\ssk

 \noi{\bf Step 4.} In view of the extra term in (\ref{H1P}),
 instead of (\ref{F11}), we will have
  \beq
  \label{F11P}
  g=u_t + g_1(u) + h(t) W'(u).
   \eeq
Therefore, by (\ref{h1}), instead of (\ref{w31}), we obtain 
holds
 \beq
 \label{w31P}
 \|u(t) - \psi_{\mu(t)}\|_{H^2} \le C_4(\|u_t\|+ h(t)).
  \eeq

\ssk

\noi{\bf Step 5.} It follows from (\ref{V11NP}) that the
convergence (\ref{V12})  takes place but not necessarily  from
below, which is not important. Then, instead of (\ref{V15}), we
have
 \beq
\label{V15P}
 \mbox{$ \frac {{\mathrm d}}{{\mathrm d}t} \, [B-
V(u(t))] = -\|u_t\|^2 + \langle h(t) W'(u),u_t \rangle
 \le - C_7[B-
V(u(t))] + C_7 h(t). $}
 \eeq
  Integrating this inequality
by using the slow decay hypotheses (\ref{h111}), (\ref{h12}), we
 conclude that, instead of the exponential estimate
(\ref{V65}), on large bounded intervals $t \in [t_k, t_k+T]$ as $k \to
\infty$, with $T \gg 1$ and $t \gg 1$,
 \beq
\label{V65P} B- V(u(t)) \le 2 h(t_k+t).
 \eeq

\ssk

\noi{\bf Step 6} needs a major revision, since Zelenyak's Lemma
\ref{Pr.Z4} applies only to exponential decay estimates.

\begin{lemma}
\label{Pr.Z4P} Assume  $(\ref{h1})$ and $(\ref{h12N})$ hold and,
for all $t >0$, 
 \beq
\label{V16P} \mbox{$
 \int\limits_t^\infty \|u_t(s)\|^2 \, {\mathrm d}s \le h(t).
 $}
 \eeq
Then
 \beq
 \mbox{$
\label{V17P} \|u(t)-u(\t)\| \le \int\limits_{t-1}^\infty \sqrt {h(s)}\, {\mathrm d}s \quad \mbox{for
all} \,\,\, 0 < t \le \t \le \infty.
 $}
 \eeq
\end{lemma}

 \noi {\em Proof.} The result is obvious if
 $|t-\t| \le 1$, where, by the H\"older inequality,
 $$
\begin{matrix}
\|u(t)-u(\t)\|
 \le \sqrt{\t-t} \,\,  \sqrt{\int\limits_t^\t \|u_t(s)\|^2 \, {\mathrm d}s} \le \sqrt{\t-t}\, \sqrt{h(t)}.
 \end{matrix}
 $$
For large intervals $\t-t>1$, we perform a partition with a
sequence of time-steps $\{\D_i, i=0,1,..., K\}$, so that,
instead of (\ref{uu1}), we obtain
 \beq
 \label{uu1P}
\begin{matrix}
\|u(t)-u(\t)\| \le \sum\limits_{i=0}^{K} \,\sqrt{ \D_i}\,\, \sqrt{
\int\limits_{t_i}^{\infty} \|u_t(s)\|^2 \, {\mathrm d}s}
= \sum\limits_{i=0}^{K}\, \sqrt{\D_i}\,\, \sqrt{h(t_i)},
 \end{matrix}
  \eeq
where $t_{i+1}= t_i + \D_i$.
Setting $\D_i = 1$ for $i=0,1,...,K-1$ and $\D_K \le 1$, by the integral test of convergence of series, we obtain
(\ref{V17P}).
 $\qed$

\ssk



\noi{\bf Step 7} remains the same, where we replace ``exponentially small"
in (a) and (b) by ``$O(\sqrt{h(t)})$-small".

\ssk

\noindent{\bf Acknowledgement. }
The second and the third author thank Department of Mathematical
Sciences of the University of Bath for its hospitality during
their visits sponsored by the INTAS Projects CERN-INTAS00-0136 and
INTAS 03-51-5007.



\ssk

\begin{small}

\noi\underline{\sc Comment}. Almost two years after publishing the
present paper in\footnote{It
is worth mentioning that this research was essentially complieted in 2005, but publishing the paper \cite{GPSZel} took quite a while,
since the first its submission to the J. Differ. Equat., after a rather 
long time,  was rejected with the critics from a Referee saying, loosely speaking,
that almost all the obtained results directly follow from Hale and
Raugel's ones \cite{Hale92} of 1992, and that the exponential
convergence in Theorem \ref{Th.C1}(iii) (not obtained elsewhere at
that time) can be also somehow easily proved. The latter is wrong,
since the exponential convergence cannot be in principle derived
from any kind of invariant manifold theory (in particular, as we
have mentioned, this has nothing to do with the classic {\em
exponential} stability of centre manifolds).
This exponential convergence is one of the main achievements of
Zelenyak's approach, which he developed for 1D second-order
parabolic equations \cite{Zel} that we here extend to more general
higher-order parabolic flows in $\ren$ and to some autonomous or
perturbed ODEs in Hilbert spaces. Note that the exponential
convergence in \cite[Theorem~1.1(iii)]{HJ07} (a Hilbert space
framework) is also proved by Zelenyak's Lemma and his related
techniques.}
$$
\mbox{Sbornik: Math., {\bf 198}:6 (2007), 817--838; see \cite{GPSZel}},
 $$
 in January 2009, the first author found that similar
Zelenyak's ideas and techniques were applied in \cite[\S~5]{HJ07}
to convergence in 1D wave equations. Actually, Zelenyak's Lemma
\ref{Pr.Z4} were introduced in \cite{HJ07} earlier, in Section 2,
and was a key ingredient of the proof of Theorem 1.1 on
exponential convergence in a Hilbert space setting (similar to our
Theorem 1.1 in an analogous Hilbert space framework explained in
Section \ref{SectHi}). It is interesting and truly remarkable that
such a growing interest to Zelenyak's fundamental idea on proving exponential convergence
 from the
1960s \cite{Zel}, in the 21st century,  occurred approximately in
the same year, about 2007, according to the publication dates.

\end{small}



\begin{thebibliography}{111}




\bibitem
 {Kond1}
 L.A. Bagirov and V.A. Kondratiev, {\em On asymptotic properties
 of solutions of diffusion equations}, Proc. Petrovskii's Seminar,
  {\bf 22} (2002), 37-70; English transl. in J.~Math. Sci. (N.Y.), {\bf 114}, No. 4 (2003), 1407--1428.




\bibitem 
 {Berger}
  M.~Berger, {\rm Nonlinearity and Functional Analysis}, Acad.
  Press, New York, 1977.



\bibitem 
{BS}
 M.S. Birman and M.Z. Solomjak, { Spectral Theory of Self-Adjoint Operators
 in Hilbert Space}, D. Reidel, Dordrecht/Tokyo, 1987.


\bibitem
 {Busca02}
 J.~Busca, M.A.~Jendoubi, and P.~Pol\'a\u{c}ik, {\em Convergence to
 equilibrium for semilinear parabolic problems}, Commun. Part. Differ. Equat.,
  {\bf 27} (2002),
 1793-1814.


\bibitem
{Chill} R. Chill, {\em On the \L{ojasiewicz}--Simon gradient
inequality}, J. Funct. Anal., {\bf 201} (2003), 572-601.


\bibitem
 {Clark}
 D.C. Clark, {\em A variant of Lusternik-Schnirelman theory},
 Indiana Univ. Math. J., {\bf 22} (1972), 65-74.


\bibitem
 {Deiml} K. Deimling, {\rm Ninlinear Functional Analysis},
 Springer-Verlag, Berlin/Tokyo, 1985.

 \bibitem
 {Sim2}
J.~Escher and  G.~Simonett, {\em A centre manifold analysis for
the Mullins--Sekerka model}, J.~Differ. Equat., {\bf 143}
 (1998), 267-292.


\bibitem 
{Fer} E.~Fereisl, F. Issard-Roch, and  H. Petzeltova, {\em A
non-smooth version of the Lojasiewicz--Simon theorem with
applications to non-local phase-field systems},
 {J. Differ. Equat.,}
{\bf 199} (2004), 1-21.



\bibitem
{Fr} A. Friedman, {Partial {D}ifferential {E}quations}, Robert E.
Krieger Publ.
  Comp., Malabar, 1983.


\bibitem
 {GalC}
 V.A. Galaktionov,
{\em Evolution completeness of separable solutions of nonlinear
diffusion equations in bounded domains}, Math. Meth. Appl. Sci.,
{\bf 27} (2004), 1755-1770.

\bibitem
 {GPSZel}
V.A.~Galaktionov, S.I.~Pohozaev, and A.E.~Shishkov,
{\em On convergence in gradient
 systems with branching of 
 equilibria}, Sbornik: Math., {\bf 198} (2007), 817--838.


 \bibitem 
{AMGV}
 V.A.~Galaktionov and J.L.~Vazquez, {A Stability Technique
  for Evolution Partial Differential Equations.
 A Dynamical Systems Approach},
Birkh\"auser, Boston/Berlin, 2004.



  \bibitem {Ha}
J.K.~Hale, {Asymptotic Behavior of Dissipative Systems}, AMS,
Providence, RI, 1988.

\bibitem
 {Hale92}
 J.K.~Hale and G.~Raugel, {\em Convergence in gradient-like systems with applications
 to PDE}, Z. angew. Math. Phys., {\bf 43} (1992), 63-124.

\bibitem
 {Haraux}
  A. Haraux and M.A. Jendoubi, {\em Decay estimates
to equilibrium for some evolution equations with an analytic
nonlinearity}, Asympt. Anal., {\bf 26} (2001), 21-36.


\bibitem
 {HJ07}
  A. Haraux and M.A. Jendoubi, {\em On the convergence of
 global and bounded solutions of some evolution equations},
 J.~Evol. Equat., {\bf 7} (2007), 449--470.



\bibitem{He}
D. Henry, {\rm Geometric {T}heory of {S}emilinear {P}arabolic
{E}quations},
  Lecture Notes in Math., Vol. {\bf 840}, Springer-Verlag, New York, 1981.


\bibitem
 {Kiel}
 H. Kielh\"ofer, {\rm Bifurcation Theory. An Introduction with Applications to
PDEs}, Springer-Verlag, New York, 2004.


\bibitem
{KrasZ} M.A. Krasnosel'skii and P.P. Zabreiko, {\rm Geometrical
Methods of Nonlinear Analysis}, Springer-Verlag, Berlin/Tokyo,
1984.


\bibitem{Lun}
A. Lunardi, {\rm Analytic {S}emigroups and {O}ptimal {R}egularity
in
  {P}arabolic {P}roblems}, Birkh\"auser, Basel/Berlin, 1995.




\bibitem
{Lus0}
L.A.~Lusternik, {\em A class of non-linear operators in Hilbert space},
Izv. Akad. Nauk. SSSR, Ser.
Mat., {\bf 3} (1938), 257-264.

\bibitem
{LusS} L. Lusternik and L. Schnirelman, {\em Sur le probl\`eme de
trois g\'eod\'esiques ferm\'ees sur les surfaces de genre O},
Comptes Rendus Acad. Sci. Paris, {\bf 189} (1929), 269-271.

\bibitem
{LusS1} L. Lusternik and L. Schnirelman, M\'ethodes topologiques
dans le probl\`emes variationels. I. Espaces \`a un nombre fini de dimensions, Hermann, Paris, 1934; Russian original: Moscow State Univ., Moscow, 1930.

\bibitem 
{Ma0}  H.~Matano, {\em Convergence of solutions of one-dimensional
semilinear parabolic equations}, {J. Math. Kyoto Univ. (JMKYAZ),}
{\bf 18} (1978), 221-227.

\bibitem
{Pazy}
 A.~Pazy, {\rm Semigroups of Linear Operators and Applications to
 Partial Differential Equations}, 2nd Edition, Springer-Verlag,
 New York, 1992.



\bibitem
 {Palis} J. Palis and W. de Melo, {\rm Geometric Theory of
 Dynamical Systems}, Springer, New York, 1982.


\bibitem
{Poh0} S.I. Pohozaev, \rm {\em On an approach to nonlinear
equations}, Soviet Math. Dokl.,
 {\bf 20} (1979), 912-916.


\bibitem
{PohFM} S.I. Pohozaev, \rm {\em The fibering method in nonlinear
variational problems}, Pitman Research Notes in Math., Vol. {\bf
365}, Pitman, 1997, pp. 35-88.

\bibitem
{Pol02} P. Pol\'a${\rm\breve{c}}$ik and F. Simondon, {\em
Nonconvergent bounded solutions of semilinear heat equations on
arbitrary domains}, J. Differ. Equat., {\bf 186} (2002), 586-610.





\bibitem
 {Sell}
 G.R. Sell and Y. You, {\rm Dynamics of Evolutionary Equations}, Springer-Verlag,
New York, 2002.


\bibitem
 {Sim1}
 G.~Simonett, {\em Centre manifolds for quasilinear
 reaction-diffusion systems}, Differ. Integr. Equat., {\bf 8}
 (1995), 753-796.

\bibitem 
{VainbergTr} M.A. Vainberg and V.A. Trenogin, {\rm Theory of
Branching of Solutions of Non-Linear Equations}, Noordhoff Int.
Publ., Leiden, 1974.



 \bibitem 
{Zel} T.I.~Zelenyak, {\em Stabilization of solutions of boundary
value problems for a second order parabolic equation with one
space variable}, {Differ. Equat.,} {\bf 4} (1968), 17-22.

\end{thebibliography}
\end{document}